\documentclass[review]{elsarticle}
\makeatletter
\def\ps@pprintTitle{%
	\let\@oddhead\@empty
	\let\@evenhead\@empty
	\def\@oddfoot{\centerline{\thepage}}%
	\let\@evenfoot\@oddfoot}
\makeatother

\usepackage{lineno,hyperref}
\usepackage{threeparttable}
\usepackage{pdflscape}
\usepackage{float}
\usepackage{booktabs}
\usepackage{amsmath}
\usepackage{amsfonts}
\usepackage{amssymb}
\usepackage{amsmath}
\usepackage{comment}
\usepackage{color}
\newtheorem{remark}{Remark}

\newcommand{\comm}[1]{}
\modulolinenumbers[5]
\journal{ }

\begin{document}

	\begin{frontmatter}
		
		\title{A comment on the article "The Schwarz alternating method in solid mechanics" by Alejandro Mota, Irina Tezaur and Coleman Alleman [Comput. Methods Appl. Mech. Engrg. 319 (2017) 19–51]}
		
		\author{Marcin Ma\'zdziarz\corref{mycorrespondingauthor}}
		\ead{mmazdz@ippt.pan.pl}
		
		\address{Institute of Fundamental Technological Research Polish Academy of Sciences,	Warsaw, Poland}
		\cortext[mycorrespondingauthor]{Corresponding author}

		
	\end{frontmatter}
	
%	\linenumbers
	
	%\section{The Elsevier article class}
	
	Recently, Alejandro Mota and Irina Tezaur and Coleman Alleman \cite{Mota201719} extended the well known Schwarz alternating method from linear to finite-deformation solid mechanics. They developed and introduced four variants of the Schwarz alternating method, presented proof of geometric convergence of the method and prepared parallel implementation applied to some examples. 
	Unfortunately, the work contains serious errors, both from the point of view of finite-deformation solid mechanics as well as mathematical elasticity.
	
	In the proof of the convergence of the Schwarz alternating method presented in Ref.\cite[Eq.36, P.3]{Mota201719} it is assumed that $\varPhi$[$\varphi$] is strictly convex with respect to the deformation gradient \textbf{\textit{F}}. Such seemingly innocuous assumption must be unfortunately ruled out:
	\begin{itemize}
		\item This assumption is unacceptable physically, enforces the uniqueness of the solution and makes description of buckling impossible, see Ref.\cite[Ch.0]{Ball1976} and references there, 
		\item moreover, it is incompatible with the property: $\varPhi$[$\varphi$]$\rightarrow$+$\infty$ as det\textbf{\textit{F}}$\rightarrow$0$^+$,
		\item as well as from objectivity condition, i.e. \textit{axiom of material frame invariance}, see Ref.\cite[Thm.4.8-1]{Ciarlet1998}.
	\end{itemize}
	Therefore, the \textit{convexity} of the energy function $\varPhi$[$\varphi$] should be replaced by the weaker \textit{polyconvexity}, i.e. sufficient condition for \textit{quasiconvexity}, see Ref.\cite{Ciarlet1998}.  
	
	The stored energy function used in numerical example in Ref.\cite[Eq.46]{Mota201719} is not strictly convex with respect to the deformation gradient \textbf{\textit{F}}, see Rmk.\ref{rem:1}, in fact it is polyconvex, and does not meet the assumptions of the Ref.\cite[Eq.36, P.3]{Mota201719}.    
	
	\begin{remark}(\textit{nonconvexity} of the energy function)
		\label{rem:1}
		
		We say that W(\textbf{F}) is convex for $\lambda$$\in$(0,1) if
		\begin{eqnarray}
		\textit{W}({\lambda}\textbf{F}_1+(1-\lambda)\textbf{F}_2)< {\lambda}\textit{W}(\textbf{F}_1)+({1-\lambda})\textit{W}(\textbf{F}_2).
		\label{eqn:C}
		\end{eqnarray}
		
		Choose \textbf{F$_1$, F$_2$, F$_{\lambda}$}$\in$$\mathbb{M}^{3\times3}$
		\begin{eqnarray}
		\textbf{F}_{1}=\left[
		\begin{array}{ccc}
		1&0&0\\
		0&1&0\\
		0&0&1
		\end{array}\right],  
		\textbf{F}_{2}=\left[
		\begin{array}{ccc}
		-1&0&0\\
		0&-1&0\\
		0&0&1
		\end{array}\right], 
		\textbf{F}_{\lambda}=\lambda\textbf{F}_{1}+(1-\lambda)\textbf{F}_{2},
		\end{eqnarray}
		
		this implies J(\textbf{F$_1$})=J(\textbf{F$_2$})=1, tr(\textbf{F$_1$$^T$F$_1$})=tr(\textbf{F$_2$$^T$F$_2$})=3,\\ J(\textbf{F}$_{\lambda}$)=(2$\lambda$-1)$^2$, tr(\textbf{F$_{\lambda}$$^T$F$_{\lambda}$})=2$\times$(2$\lambda$-1)$^2$+1.
		
		Stored energy function in Ref.\cite[Eq.46]{Mota201719}:
		\begin{eqnarray}
		\textit{W}(\textbf{C})=\frac{\kappa}{4}(J^2-2logJ-1)+\frac{\mu}{2}(J^{-2/3}tr\textbf{C}-3), where~J=det(\textbf{F}), \textbf{C}=\textbf{F}^T\textbf{F}.
		\label{eqn:NH}
		\end{eqnarray}
		
		For convenience choose {$\kappa$}=4, {$\mu$}=2. Omitting some cumbersome calculations for {$\lambda$}$\in$(0,1)$\backslash$1/2 
		we get \textit{W}(\textbf{F}$_{1}$)=\textit{W}(\textbf{F}$_{2}$)=0,\textit{W}(\textbf{F}$_{\lambda}$)$>$0 and contradiction of convexity condition in Eq.\ref{eqn:C}.
	\end{remark}
	
	Let's assume now that we replace the \textit{convexity} assumption in the proof of the convergence of the Schwarz alternating method in Ref.\cite[Eq.36, P.3]{Mota201719} by the more suitable  weaker \textit{polyconvexity} and try to apply the proposed Schwarz alternating method to a problem of Buckling of a Rod analysed in Ref.\cite[Ch.9]{Ball1976}. It is shown \textit{ibidem} that, "\textit{sufficiently long rods of arbitrary cross-section will exhibit nonuniqueness}" and buckle. If we divide such a sufficiently long rod into shorten domains, as it is done in the proposed Schwarz alternating method, this condition of sufficient length will be not fulfilled and we will not get buckling. Thus, we will receive a solution other than for one domain. Consequently the Schwarz alternating method will give a completely different qualitative solution, that rather excludes it from the use in application of finite-deformation solid mechanics where nonuniqueness can usually occur. 
	
	\section*{References}
	
	\bibliography{CMAMEbibfile}
	
\end{document}